\documentclass[12pt]{article}

\usepackage{tikz}

\usetikzlibrary{calc,through,backgrounds}
\usetikzlibrary{decorations.pathreplacing}

\usepackage{multicol}

\usepackage{bm,amssymb,amsthm,amsfonts,enumerate,amsmath,latexsym,cite,color, psfrag,graphicx,ifpdf,pgfkeys,pgfopts,xcolor,extarrows}

\newtheorem{theo}{Theorem}[section]

\newtheorem{lem} [theo]{Lemma}
\newtheorem{coro}[theo]{Corollary}
\newtheorem{prop}[theo]{Proposition}

\newtheorem{re}[theo]{Remark}

\allowdisplaybreaks

\makeatletter \@addtoreset{equation}{section}
\@addtoreset{theo}{}\makeatother

\usepackage{hyperref}

\def\qed{\hfill \rule{4pt}{7pt}}
\def\pf{\noindent {\it Proof.} }

\begin{document}
\begin{center}
{\Large\bf An Involution on Semistandard  Skyline Fillings}

\vskip 4mm
{\small  NEIL J.Y. FAN, PETER L. GUO, NICOLAS Y. LIU}

\end{center}

\noindent{A{\scriptsize BSTRACT}.}
Non-attacking skyline fillings     were used  by
Haglund,  Haiman  and  Loehr to establish a combinatorial formula for nonsymmetric
Macdonald polynomials. Semistandard skyline fillings are non-attacking
skyline fillings with both major index and coinversion number   equal to zero, which
 serve as a combinatorial model for   key polynomials. In this paper, we
construct  an involution on semistandard skyline fillings. This involution can be
viewed as a vast generalization of the classical  Bender--Knuth involution.
 As an application, we obtain that
 semistandard skyline fillings are compatible with the
Demazure operators, offering a new combinatorial  proof that
nonsymmetric Macdonald polynomials specialize to key polynomials.

%
%

\section{Introduction}

The key polynomials  $\kappa_\alpha(x)$ associated to compositions $\alpha\in \mathbb{Z}_{\geq 0}^n$,
also called   Demazure characters,
 are characters of the Demazure
modules for the general linear groups \cite{Dem-1,Dem-2}.
These polynomials are  defined based on the {\it Demazure operator} $\pi_i=\partial_ix_i$.
Here, 
$\partial_i$ is the {\it divided difference operator} sending
  a polynomial $f(x)\in \mathbb{Z}[x_1,x_2,\ldots,x_n]$
  to
  \[\partial_i(f(x))=\frac{f(x)    -s_i f(x)}{x_i-x_{i+1}},\]
where $s_i f(x)$ is obtained from $f(x)$ by interchanging $x_i$ and $x_{i+1}$. Precisely, if $\alpha=(\alpha_1,\alpha_2,\ldots, \alpha_n)$ is a partition (i.e., $\alpha_1\geq \alpha_2\geq \cdots\geq \alpha_n$), then set
$\kappa_\alpha(x)=x_1^{\alpha_1}x_2^{\alpha_2} \cdots x_n^{\alpha_n}.$
Otherwise, choose an index  $i$ such that $\alpha_i<\alpha_{i+1}$, and let $\alpha'$ be
 obtained from $\alpha$ by interchanging $\alpha_i$ and $\alpha_{i+1}$. Set
\begin{equation}\label{induc}
\kappa_\alpha(x)=\pi_i (\kappa_{\alpha'}(x))  =\partial_i(x_i\kappa_{\alpha'}(x)).
\end{equation}
The above definition is independent of the choice of
  $i$  since the Demazure operators  satisfy the
Coxeter relations: $\pi_i \pi_j=\pi_j \pi_i$ for $|i-j|>1$, and
$\pi_i \pi_{i+1} \pi_i= \pi_{i+1}\pi_i  \pi_{i+1}$.

Combinatorial constructions of key polynomials have received extensive attention.
Lascoux and   Sch\"utzenberger  \cite{Las-1} initiated a combinatorial treatment
for key polynomials, showing that
$\kappa_\alpha(x)$ can be expressed  as a sum of the weights of {\it semistandard Young tableaux}  with right key
bounded by a certain key, see also Reiner and  Shimozono \cite{Rei}. Kohnert \cite{Koh} provided
 a  diagram interpretation, now called {\it Kohnert diagrams},  for
$\kappa_\alpha(x)$ by giving a bijection between Kohnert diagrams
   and semistandard Young tableaux with bounded right key, see also  Assaf and  Quijada
\cite{Ass-1}.  Assaf and Searles \cite{Ass-2} defined  the notion of
{\it Kohnert tableaux}, and established a bijection between Kohnert tableaux and
  Kohnert diagrams.
Assaf \cite{Ass} introduced    {\it semistandard key tableaux}
to give another combinatorial construction of
 key polynomials.
Comparing the definitions of Kohnert tableaux \cite{Ass-2}
 and semistandard key tableaux \cite{Ass}, it can be  checked  that
  Kohnert tableaux  are in bijection with
 semistandard key tableaux \cite{X}.

 On the other hand, Ion \cite{Ion} proved that $\kappa_\alpha(x)$  can be realized
 as a specialization of the nonsymmetrict Macdonald polynomial $E_\alpha(x; q, t)$
at $q=t=0$, namely,
\begin{align}\label{keye}
\kappa_\alpha(x)=E_\alpha(x; q=0, t=0).
\end{align}
 Nonsymmetrict Macdonald polynomials   were
  developed by  Opdam \cite{Opd}, Macdonald
\cite{Mac} and  Cherednik \cite{Che}.
Haglund,  Haiman  and  Loehr \cite{Hag} established a combinatorial formula
for $E_\alpha(x; q, t)$ in terms of {\it non-attacking skyline fillings} for  $\alpha$
 with $x$ counting the weights,
$q$ counting the coinversion numbers, and  $t$ counting the  major indices.
This, combined with \eqref{keye}, implies that $\kappa_\alpha(x)$ is a weighted counting
of non-attacking skyline fillings for $\alpha$  with both major index and
coinversion number equal to zero. Such non-attacking fillings
are called {\it semistandard skyline fillings}, see for example
 Monical \cite{Mon-2} or Monical,  Pechenik and  Searles \cite{Mon}.

As pointed out by Assaf \cite[Proposition 3.1]{Ass},
semistandard skyline fillings are exactly
semistandard key tableaux.  Assaf  \cite{Ass} also
gave an alternative proof of the relation
\eqref{keye}
 based on
an expansion of   key polynomials into
 fundamental slide polynomials
established in \cite[Corollary 3.16]{Ass-x}.

In  Figure \ref{Fig-skyline}, we summarize   the relationships  among  the above mentioned
combinatorial constructions  of key polynomials.

\begin{figure}[h]
\begin{center}
\begin{tikzpicture}[scale=0.5]

 \node at (-10mm,0mm) [align=center,draw,text width=4.1cm]
{semistandard Young tableaux with bounded right key};

\node at (58mm,1mm)
 {$\underleftrightarrow{\hspace{2.2cm}}$};
\node at (58mm,8mm) {\footnotesize bijection \cite{Koh}};

 \node at (110mm,0mm)  [align=center,draw,text width=2.5cm]
{Kohnert diagrams};

\node at (160mm,1mm)
 {$\underleftrightarrow{\hspace{2.1cm}}$};
\node at (160mm,8mm) {\footnotesize bijection \cite{Ass-2}};

\node at (210mm,0mm)  [align=center,draw,text width=2.5cm]
{Kohnert tableaux };

\node at (187mm,-30mm){\footnotesize bijection  \cite{X}};

\draw[<->] (210mm,-45mm) -- (210mm,-15mm);

\node at (200mm,-60mm)  [align=center,draw,text width=3.1cm]
{semistandard key tableaux};

\node at (142mm,-58mm)
 {$\underleftrightarrow{\hspace{2.3cm}}$};
\node at (142mm,-51mm) {\footnotesize the same \cite{Ass}};

\node at (80mm,-60mm)  [align=center,draw,text width=3.5cm]
{semistandard skyline fillings};

\end{tikzpicture}
\end{center}
\vspace{-2mm}
\caption{Relations among some combinatorial models of key polynomials. }\label{Fig-skyline}
\end{figure}
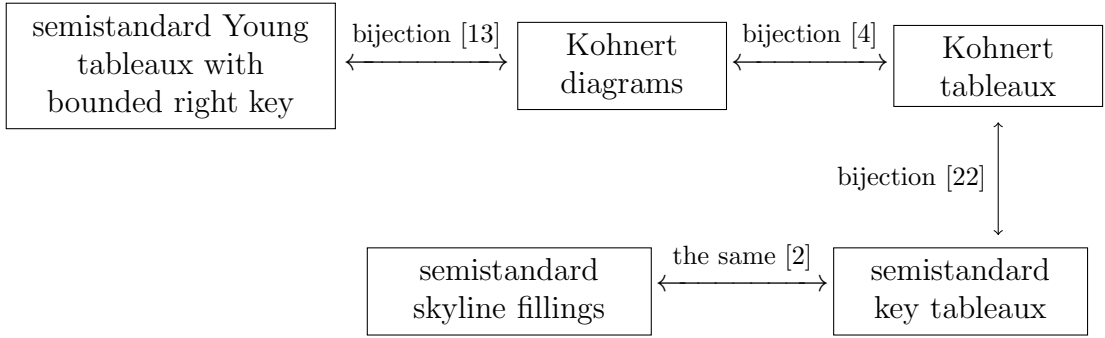

In this paper, we construct   an involution on   semistandard skyline fillings.
This involution significantly generalizes the classical Bender--Knuth involution
from semistandard Young tableaux to semistandard skyline fillings.
To be specific, when the parts of $\alpha$ are weakly increasing, the key polynomial
$\kappa_{\alpha}(x)$ becomes  a Schur polynomial, and our involution
specifies to the   Bender--Knuth involution \cite{Ben}.
As an application, it immediately follows  that  semistandard skyline
fillings are compatible with the Demazure operators. This
provides a new simple   proof of the relation \eqref{keye} between
nonsymmetric Macdonald polynomials and key polynomials.

We end this section with some remarks,
which  are also parts of the motivation of this paper.
The first two remarks concern
combinatorial realizations of the Demazure operators in two
formulas for
key polynomials.
As comparison, our treatment avoids using the properties of
the ``plactic congruence'' or the RSK insertion algorithm.

\vspace{5pt}

\begin{re}\label{Re-1}
As aforementioned, Lascoux and   Sch\"utzenberger  \cite{Las-1}
established the following formula of $\kappa_\alpha(x)$:
\begin{equation}\label{MMm}
\kappa_\alpha(x)=\sum_{T\in \mathcal{T}(\alpha)}x^{T},
\end{equation}
where
\[\mathcal{T}(\alpha)=\{T\colon  K_{+}(T)\,\leq\, \mathrm{key}
(\alpha)\}\]
is the set of  semistandard Young tableaux $T$ of shape $\lambda(\alpha)$
(namely, the partition by rearranging the
parts of $\alpha$)  such that
the right key $K_{+}(T)$ of $T$ is entry-wisely less than or equal to   $\mathrm{key}(\alpha)$. Here,    $\mathrm{key}(\alpha)$ is the semistandard  Young tableau of shape $\lambda(\alpha)$ whose first $\alpha_j$ columns  contain the letter $j$ for all $j$.
The definition
of the right key $K_{+}(T)$ of  $T$ is quite subtle, and  we
 refer the reader to \cite{Rei} for a detailed description.
Lascoux and   Sch\"utzenberger  \cite{Las-1} proved  \eqref{MMm}  by introducing a combinatorial version
of the Demazure operator on   semistandard Young
tableaux,
see also \cite[Section 4]{Len} or \cite[Section 2]{Mas}.
The proof relies heavily on properties of   ``plactic congruence''.

\end{re}

\vspace{5pt}

\begin{re}\label{Re-2}
Another combinatorial description of $\kappa_\alpha(x)$ also due to
 Lascoux and   Sch\"utzenberger  \cite{Las-1}  is as
follows:
\begin{equation}\label{LS-2}
\kappa_\alpha(x)=\sum_{u\in \mathcal{W}(\alpha)}x^{u},
\end{equation}
where $\mathcal{W}(\alpha)$ is the set of certain flagged  row-frank words, see   \cite{Rei}
for the precise   definition. Reiner and Shimozono
\cite{Rei}  gave a  proof  of \eqref{LS-2} by showing combinatorially that the set
$\mathcal{W}(\alpha)$ is compatible with  the Demazure operators, namely,
\begin{equation}\label{RS-P}
\pi_r\left(
\sum_{v\in \mathcal{W}(\alpha')}x^v\right)=\sum_{u\in \mathcal{W}(\alpha)}x^u,
\end{equation}
where $\alpha'$ is a composition
 obtained from $\alpha$ by interchanging $\alpha_i$ and $\alpha_{i+1}$ with
 $\alpha_i<\alpha_{i+1}$.  The proof of \eqref{RS-P} uses  implicit
 properties of the ``$r$-paring'' operation    and the RSK
  insertion algorithm.
\end{re}

\vspace{5pt}

 The third remark discusses whether the involution
 could be extended to give a combinatorial proof of a
 formula for Lascoux polynomials.

\begin{re}
The Lascoux polynomial $L_\alpha(x)$ is the $K$-theory analog of the key
polynomial $\kappa_\alpha(x)$. In fact, $\kappa_\alpha(x)$ can be
obtained from $L_\alpha(x)$ by extracting the lowest degree component.
Monical \cite{Mon-2} gave a conjectural
 formula for
$L_\alpha(x)$ in terms of semistandard set-valued
skyline fillings for $\alpha$, see also \cite{Mon}. This conjecture was recently
confirmed by Buciumas,  Scrimshaw  and   Weber
\cite{Buc}
using the colored five-vertex model. We do not know  if our
involution could be extended  to a set-valued version  to give a
 combinatorial proof of Monical's formula for Lascoux polynomials.
We point out that   Miller \cite{Mil} showed
that the  classical reduced pipe dreams (or, RC-graphs)   for
the constructions of
 Schubert polynomials   are compatible with
the divided difference operators by using  the mitosis algorithm \cite{KnMi}.
Grothendieck polynomials are $K$-theory analogs of Schubert polynomials.
The mitosis algorithm   has been extended by Tyurin \cite{Tyu}
to   non-reduced pipe dreams for the construction of   Grothendieck polynomials.
\end{re}

This paper is structured as follows. In Section \ref{sec-2}, we
construct the promised  involution on simistandard skyline fillings.
The key point in the construction is a classification of two consecutive
entries  appearing in a simistandard skyline filling.
In Section \ref{sec-3},  we apply this involution to give a
combinatorial explanation that
semistandard skyline fillings are compatible with the
Demazure operators.

\section{The involution on simistandard skyline fillings}\label{sec-2}

In this section, we aim to construct an involution   on semistandard skyline fillings.
Let us  start  with an overview of
semistandard skyline fillings.

\subsection{Semistandard skyline fillings}

Let $\alpha=(\alpha_1,\alpha_2,\ldots,\alpha_n)\in \mathbb{Z}_{\geq  0}^n$ be a  composition.
The skyline diagram $D(\alpha)$ of $\alpha$ is  a left-justified array
with $\alpha_i$ boxes in row $i$. Here, the row indices increase from
top to bottom, and the column indices increase from left to right.
We use $(i,j)$ to denote the box in row $i$ and column $j$.
For example,  Figure \ref{Fig-skyline} is an illustration of the skyline diagram of  $(1,3,0,2)$.

\vspace{2mm}
\begin{figure}[h]
\begin{center}
\begin{tikzpicture}[scale=0.5]
      \draw [-, shift={(0,0)}] (0,0)--(0,4);
      \filldraw [fill=gray, fill opacity=0, shift={(0,0)}] (0,0) rectangle (1,1);
      \filldraw [fill=gray, fill opacity=0, shift={(0,0)}] (1,0) rectangle (2,1);
      \filldraw [fill=gray, fill opacity=0, shift={(0,0)}] (0,2) rectangle (1,3);
      \filldraw [fill=gray, fill opacity=0, shift={(0,0)}] (1,2) rectangle (2,3);
      \filldraw [fill=gray, fill opacity=0, shift={(0,0)}] (2,2) rectangle (3,3);
      \filldraw [fill=gray, fill opacity=0, shift={(0,0)}] (0,3) rectangle (1,4);
\end{tikzpicture}
\end{center}
\vspace{-5mm}
\caption{The skyline diagram   for $\alpha=(1,3,0,2)$. }\label{Fig-skyline}
\end{figure}
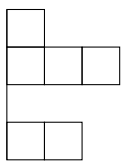

A  {\it skyline filling} for $\alpha$ is an assignment $F$ of positive integers into
the boxes of $D(\alpha)$. Write $F(i,j)$ for the entry filled in the box $(i,j)$.
A skyline  filling $F$ is {\it non-attacking} if the entries are distinct in
each column, and  $F(i, j)\neq F(i',j+1)$ for two boxes in consecutive columns with
$i>i'$.
Non-attacking  skyline fillings
were used by Haglund,  Haiman  and  Loehr \cite{Hag}   to
give a combinatorial formula for the nonsymmtric Macdonald polynomial $E_\alpha(x; q, t)$.
A {\it semistandard skyline filling}   is
a non-attaking skyline filling with both
major index and coinversion number equal to zero.

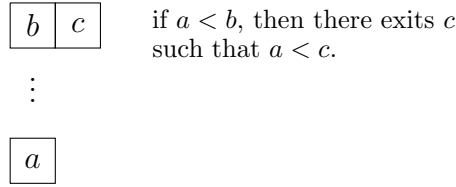
\begin{figure}[h]
\begin{center}
\begin{tikzpicture}[scale=0.6]
\draw [-] (0,3)--(1,3)--(1,2)--(0,2)--(0,3);
\draw [-] (1,3)--(2,3)--(2,2)--(1,2)--(1,3);
\draw [-] (0,-1)--(1,-1)--(1,0)--(0,0)--(0,-1);
\draw (0.5,2.5) node {$b$};
\draw (1.5,2.5) node {$c$};
\draw (0.5,-0.5) node {$a$};
\draw (0.5,1.25) node {$\vdots$};
\node at (6.5,2.6) {\footnotesize if $a<b$, then there exits $c$};
\node at (5.2,2) {\footnotesize  such that $a<c$.};
\end{tikzpicture}
\end{center}
\caption{An illustration of the condition \eqref{ifv}.}\label{Fig-iv}
\end{figure}

As mentioned in
Introduction, Assaf \cite[Proposition 3.1]{Ass} showed that
 semistandard skyline fillings are exactly
 semistandard key tableaux.
In this paper, we adopt  the equivalent definition due to Assaf, which is relatively easier to describe.
To be specific, a {  semistandard skyline filling}
is a skyline filling   such that
\begin{enumerate}[(i)]

\item \label{ifi}
the entries are weakly decreasing from left to right in each row;

\item  \label{ifii}
each entry cannot exceed its row index (flag condition);

\item  \label{if3i}
the entries are distinct in each column;

\item   \label{ifv}
if some entry $a$ is below  and in the same column
as an entry $b$ with $a<b$, then there is an entry to the immediately  right of $b$, say $c$,
such that $a<c$. See Figure \ref{Fig-iv} for an illustration.

\end{enumerate}

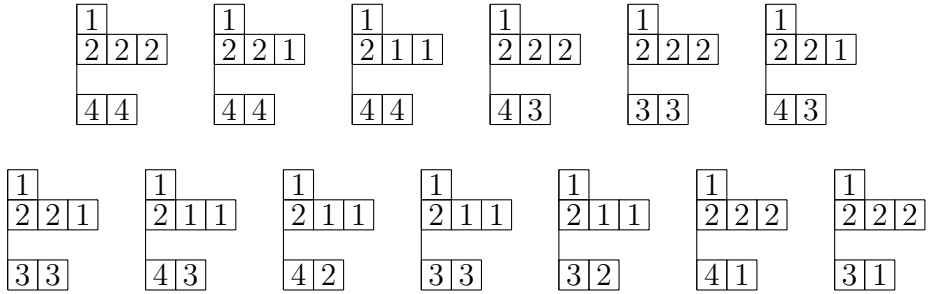
\begin{figure}[h]
\begin{center}
      \begin{tikzpicture}[scale=0.4]
      \draw (0,0) -- (0,4);
      \draw (0,0) grid (1,1); \draw (1,0) grid (2,1);
      \draw (0,2) grid (1,3); \draw (1,2) grid (2,3); \draw (2,2) grid (3,3);
      \draw (0,3) grid (1,4);
      \draw (0.5,0.5) node {$4$}; \draw (1.5,0.5) node {$4$};
      \draw (0.5,2.5) node {$2$}; \draw (1.5,2.5) node {$2$};\draw (2.5,2.5) node {$2$};
      \draw (0.5,3.5) node {$1$};
      \end{tikzpicture}
      \quad
      \begin{tikzpicture}[scale=0.4]
      \draw (0,0) -- (0,4);
      \draw (0,0) grid (1,1); \draw (1,0) grid (2,1);
      \draw (0,2) grid (1,3); \draw (1,2) grid (2,3); \draw (2,2) grid (3,3);
      \draw (0,3) grid (1,4);
      \draw (0.5,0.5) node {$4$}; \draw (1.5,0.5) node {$4$};
      \draw (0.5,2.5) node {$2$}; \draw (1.5,2.5) node {$2$};\draw (2.5,2.5) node {$1$};
      \draw (0.5,3.5) node {$1$};
      \end{tikzpicture}
      \quad
      \begin{tikzpicture}[scale=0.4]
      \draw (0,0) -- (0,4);
      \draw (0,0) grid (1,1); \draw (1,0) grid (2,1);
      \draw (0,2) grid (1,3); \draw (1,2) grid (2,3); \draw (2,2) grid (3,3);
      \draw (0,3) grid (1,4);
      \draw (0.5,0.5) node {$4$}; \draw (1.5,0.5) node {$4$};
      \draw (0.5,2.5) node {$2$}; \draw (1.5,2.5) node {$1$};\draw (2.5,2.5) node {$1$};
      \draw (0.5,3.5) node {$1$};
      \end{tikzpicture}
      \quad
      \begin{tikzpicture}[scale=0.4]
      \draw (0,0) -- (0,4);
      \draw (0,0) grid (1,1); \draw (1,0) grid (2,1);
      \draw (0,2) grid (1,3); \draw (1,2) grid (2,3); \draw (2,2) grid (3,3);
      \draw (0,3) grid (1,4);
      \draw (0.5,0.5) node {$4$}; \draw (1.5,0.5) node {$3$};
      \draw (0.5,2.5) node {$2$}; \draw (1.5,2.5) node {$2$};\draw (2.5,2.5) node {$2$};
      \draw (0.5,3.5) node {$1$};
      \end{tikzpicture}
      \quad
      \begin{tikzpicture}[scale=0.4]
      \draw (0,0) -- (0,4);
      \draw (0,0) grid (1,1); \draw (1,0) grid (2,1);
      \draw (0,2) grid (1,3); \draw (1,2) grid (2,3); \draw (2,2) grid (3,3);
      \draw (0,3) grid (1,4);
      \draw (0.5,0.5) node {$3$}; \draw (1.5,0.5) node {$3$};
      \draw (0.5,2.5) node {$2$}; \draw (1.5,2.5) node {$2$};\draw (2.5,2.5) node {$2$};
      \draw (0.5,3.5) node {$1$};
      \end{tikzpicture}
      \quad
      \begin{tikzpicture}[scale=0.4]
      \draw (0,0) -- (0,4);
      \draw (0,0) grid (1,1); \draw (1,0) grid (2,1);
      \draw (0,2) grid (1,3); \draw (1,2) grid (2,3); \draw (2,2) grid (3,3);
      \draw (0,3) grid (1,4);
      \draw (0.5,0.5) node {$4$}; \draw (1.5,0.5) node {$3$};
      \draw (0.5,2.5) node {$2$}; \draw (1.5,2.5) node {$2$};\draw (2.5,2.5) node {$1$};
      \draw (0.5,3.5) node {$1$};
      \end{tikzpicture}\quad\\[12pt]
      \begin{tikzpicture}[scale=0.4]
      \draw (0,0) -- (0,4);
      \draw (0,0) grid (1,1); \draw (1,0) grid (2,1);
      \draw (0,2) grid (1,3); \draw (1,2) grid (2,3); \draw (2,2) grid (3,3);
      \draw (0,3) grid (1,4);
      \draw (0.5,0.5) node {$3$}; \draw (1.5,0.5) node {$3$};
      \draw (0.5,2.5) node {$2$}; \draw (1.5,2.5) node {$2$};\draw (2.5,2.5) node {$1$};
      \draw (0.5,3.5) node {$1$};
      \end{tikzpicture}
      \quad
      \begin{tikzpicture}[scale=0.4]
      \draw (0,0) -- (0,4);
      \draw (0,0) grid (1,1); \draw (1,0) grid (2,1);
      \draw (0,2) grid (1,3); \draw (1,2) grid (2,3); \draw (2,2) grid (3,3);
      \draw (0,3) grid (1,4);
      \draw (0.5,0.5) node {$4$}; \draw (1.5,0.5) node {$3$};
      \draw (0.5,2.5) node {$2$}; \draw (1.5,2.5) node {$1$};\draw (2.5,2.5) node {$1$};
      \draw (0.5,3.5) node {$1$};
      \end{tikzpicture}
      \quad
      \begin{tikzpicture}[scale=0.4]
      \draw (0,0) -- (0,4);
      \draw (0,0) grid (1,1); \draw (1,0) grid (2,1);
      \draw (0,2) grid (1,3); \draw (1,2) grid (2,3); \draw (2,2) grid (3,3);
      \draw (0,3) grid (1,4);
      \draw (0.5,0.5) node {$4$}; \draw (1.5,0.5) node {$2$};
      \draw (0.5,2.5) node {$2$}; \draw (1.5,2.5) node {$1$};\draw (2.5,2.5) node {$1$};
      \draw (0.5,3.5) node {$1$};
      \end{tikzpicture}
      \quad
     \begin{tikzpicture}[scale=0.4]
      \draw (0,0) -- (0,4);
      \draw (0,0) grid (1,1); \draw (1,0) grid (2,1);
      \draw (0,2) grid (1,3); \draw (1,2) grid (2,3); \draw (2,2) grid (3,3);
      \draw (0,3) grid (1,4);
      \draw (0.5,0.5) node {$3$}; \draw (1.5,0.5) node {$3$};
      \draw (0.5,2.5) node {$2$}; \draw (1.5,2.5) node {$1$};\draw (2.5,2.5) node {$1$};
      \draw (0.5,3.5) node {$1$};
      \end{tikzpicture}
      \quad
      \begin{tikzpicture}[scale=0.4]
      \draw (0,0) -- (0,4);
      \draw (0,0) grid (1,1); \draw (1,0) grid (2,1);
      \draw (0,2) grid (1,3); \draw (1,2) grid (2,3); \draw (2,2) grid (3,3);
      \draw (0,3) grid (1,4);
      \draw (0.5,0.5) node {$3$}; \draw (1.5,0.5) node {$2$};
      \draw (0.5,2.5) node {$2$}; \draw (1.5,2.5) node {$1$};\draw (2.5,2.5) node {$1$};
      \draw (0.5,3.5) node {$1$};
      \end{tikzpicture}
      \quad
      \begin{tikzpicture}[scale=0.4]
      \draw (0,0) -- (0,4);
      \draw (0,0) grid (1,1); \draw (1,0) grid (2,1);
      \draw (0,2) grid (1,3); \draw (1,2) grid (2,3); \draw (2,2) grid (3,3);
      \draw (0,3) grid (1,4);
      \draw (0.5,0.5) node {$4$}; \draw (1.5,0.5) node {$1$};
      \draw (0.5,2.5) node {$2$}; \draw (1.5,2.5) node {$2$};\draw (2.5,2.5) node {$2$};
      \draw (0.5,3.5) node {$1$};
      \end{tikzpicture}
      \quad
      \begin{tikzpicture}[scale=0.4]
      \draw (0,0) -- (0,4);
      \draw (0,0) grid (1,1); \draw (1,0) grid (2,1);
      \draw (0,2) grid (1,3); \draw (1,2) grid (2,3); \draw (2,2) grid (3,3);
      \draw (0,3) grid (1,4);
      \draw (0.5,0.5) node {$3$}; \draw (1.5,0.5) node {$1$};
      \draw (0.5,2.5) node {$2$}; \draw (1.5,2.5) node {$2$};\draw (2.5,2.5) node {$2$};
      \draw (0.5,3.5) node {$1$};
      \end{tikzpicture}
\end{center}
\vspace{-2mm}
\caption{Semistandard skyline fillings for $(1,3,0,2)$.}\label{Fig-SSF}
\end{figure}

Let $\mathrm{SSF}(\alpha)$ denote the set of semistandard skyline fillings for $\alpha$.
For example, there are
 13 semistandard skyline fillings for  $\alpha=(1,3,0,2)$,
 as illustrated in Figure \ref{Fig-SSF}.
 \vspace{2mm}

For $F\in \mathrm{SSF}(\alpha)$, write $x^F$ for the monomial
generated by $F$, namely,
\[x^F=\prod_{(i,j)\in D(\alpha)}x_{F(i,j)}.\]
Using the combinatorial  formula for the nonsymmetric
Macdonald polynomial $E_\alpha(x; q,t)$ \cite{Hag} together with the relation
\eqref{keye}, it follows  that
\begin{equation}\label{key}
\kappa(x)=\sum_{F\in \mathrm{SSF}(\alpha)} x^F,
\end{equation}
see Assaf \cite{Ass} for an alternative proof of \eqref{key}.

The following lemma will be used later, which is in fact the non-attacking condition satisfied by
a semistandard skyline filling. Since we adopt the definition
of semistandard skyline fillings due to  Assaf \cite{Ass},
we  give a simple   proof here.

\begin{lem}\label{BB-f}
For $F\in \mathrm{SSF}(\alpha)$, if $F(i,j)=F(i',j+1)$, then we have $i\le i'$.
\end{lem}

\pf
Suppose otherwise that  $i>i'$.
Let $a=F(i', j)$. Since the entries in each row of  $F$ are weakly decreasing
 and  the entries in each column of $F$ are distinct,
we obtain that $a>F(i,j)$. However, by item \eqref{ifv}     in the definition of a semistandard skyline filling,
$F(i',j+1)$ would be strictly larger than $t$, leading to a contradiction.
\qed

In the rest of this section, we  always  fix a composition $\alpha\in
\mathbb{Z}_{\geq 0}^n$, a row index $r$ and a
 positive integer $t$.
Our involution,   denoted   $\Phi_{r,t}$, is parameterized by  $r$ and $t$. The construction of $\Phi_{r,t}$ is based on two operators: the lowering operator
  $L_{r,t}$ and the raising operator   $R_{r,t}$.
  To define these two operators,
we   classify the entries
$t$ and $t+1$ appearing in   $F\in \mathrm{SSF}(\alpha)$.
This classification is the main ingredient in the construction
of $\Phi_{r,t}$.

\subsection{Classification of two consecutive  entries}\label{sub-2-2}

Let $F\in \mathrm{SSF}(\alpha)$. Recall that the entries in every column of $F$ are
distinct. So every column of $F$
contains at most one $t$ and at most one $t+1$.
An entry  $t$ (respectively, $t+1$)
in $F$ is called {\it paired} if  there is also a $t+1$ (respectively, $t$) in the same column, otherwise it is called {\it unpaired}.
The unpaired entries $t$ or $t+1$ in $F$ will be  further classified into pseudo-free or free.

The pseudo-free entries appear in pairs (but not in the same column). Precisely,
two unpaired entries $t$ and $t+1$ in $F$ are   called {\it pseudo-free}
if
\begin{enumerate}[(1)]
\item the entry $t+1$ appears to the upper right of the entry $t$; and

\item \label{PP-x}
each column of $F$ between
this pair of $t$ and $t+1$ contains both $t$ and $t+1$ (By (1) and Lemma \ref{BB-f},
in each such column, the entry $t+1$ must be above the entry $t$).
\end{enumerate}
An   entry $t$ or $t+1$  is called  {\it free} if it is neither paired nor pseudo-free.

For example, Figure \ref{Fig-psf} depicts  a
 semistandard skyline filling. Choose  $t=2$. We use boldface to
 signify the entries $t$ and $t+1$, where   free entries are underlined and
 pseudo-free entries are  circled.

\begin{figure}[h]
\begin{center}
\begin{tikzpicture}[scale=0.6]
\draw (0,0) grid (3,3); \draw (3,1) grid (4,3);
\draw (4,2) grid (8,3); \draw (0,3) grid (11,4);
\draw (0,6) grid (4,7);
\draw (0,-1) grid (6,0); \draw (-1,-1) grid (0,7);
\draw (4,1) grid (10,2); \draw (0,5) grid (1,6);
\draw (-0.5,6.5) node {$1$}; \draw (-0.5,5.5) node {$\bm2$};
\draw (-0.5,4.5) node {$\bm3$}; \draw (-0.5,3.5) node {$4$};
\draw (-0.5,2.5) node {$5$}; \draw (-0.5,1.5) node {$6$};
\draw (-0.5,0.5) node {$7$}; \draw (-0.5,-0.5) node {$8$};
\draw (0.5,6.5) node {$1$}; \draw (1.5,6.5) node {$1$};
\draw (2.5,6.5) node {$1$}; \draw (3.5,6.5) node {$1$};
\draw (0.5,5.5) node {$\underline{\bm2}$};
\draw (0.5,3.5) node {$4$}; \draw (1.5,3.5) node {$4$};
\draw (2.5,3.5) node {$4$}; \draw (3.5,3.5) node {$4$};
\draw (4.5,3.5) node {$4$}; \draw (5.5,3.5) node {$4$};
\draw (6.5,3.5) node {$4$}; \draw (7.5,3.5) node {$4$};
\draw (8.5,3.5) node {$4$}; \draw (9.5,3.5) node {$\bm3$};
\draw (10.5,3.5) node {$\bm3$};
\draw (0.5,2.5) node {$5$}; \draw (1.5,2.5) node {$5$};
\draw (2.5,2.5) node {$\bm3$}; \draw (3.5,2.5) node {${\bm3}$};
\draw (4.5,2.5) node {$\underline{\bm3}$}; \draw (5.5,2.5) node {$\bm3$};
\draw (6.5,2.5) node {$\underline{\bm3}$}; \draw (7.5,2.5) node {$\underline{\bm2}$};
\draw (0.5,1.5) node {$6$}; \draw (1.5,1.5) node {$6$};
\draw (2.5,1.5) node {$5$}; \draw (3.5,1.5) node {$5$};
\draw (4.5,1.5) node {$5$}; \draw (5.5,1.5) node {$5$};
\draw (6.5,1.5) node {$5$}; \draw (7.5,1.5) node {$5$};
\draw (8.5,1.5) node {$\bm2$}; \draw (9.5,1.5) node {$\bm2$};
\draw (0.5,0.5) node {$7$}; \draw (1.5,0.5) node {$\bm2$};
\draw (2.5,0.5) node {$\bm2$};
\draw (0.5,-0.5) node {$8$}; \draw (1.5,-0.5) node {$8$};
\draw (2.5,-0.5) node {$7$}; \draw (3.5,-0.5) node {$7$};
\draw (4.5,-0.5) node {$6$}; \draw (5.5,-0.5) node {$\bm2$};
\draw (1.5,0.5) circle (0.4cm); \draw (3.5,2.5) circle (0.4cm);
\draw (10.5,3.5) circle (0.4cm); \draw (8.5,1.5) circle (0.4cm);
\end{tikzpicture}
\end{center}
\caption{Classification of the entries $t$ and $t+1$.  }\label{Fig-psf}
\end{figure}

%


The following lemma gives an alternative  characterization of pseudo-free entries.

\begin{lem}\label{lem-1}
Let $F\in \mathrm{SSF}(\alpha)$, and $F(i,j)=t$ be an unpaired entry.
Assume that $i'$ is a row index such that $i'<i$.
Then the following are equivalent:
\begin{enumerate}[(1)]
\item  $t$ is pseudo-free and its  associated  pseudo-free entry
 $t+1$ is in   row $i'$;

\item  $F(i',j+1)=t+1$.
\end{enumerate}
\end{lem}

\pf We first show that  (2) implies (1).
If the entry $t+1$  in the box $(i',j+1)$ of $F$ is unpaired, then
it is the pseudo-free entry associated with  $t$.
 Otherwise,   column $j+1$  contains an entry $t$,
 we assert that
 \begin{equation}\label{PPP}
 F(i', j+2)=t+1.
 \end{equation}
By Lemma \ref{BB-f},
 the entry $t$ in column $j+1$ lies in row $i$ or
in a row below row $i$.
It follows that in column $j+1$, the entry $t+1$ is above
the entry $t$. By item \eqref{ifv} in the definition of a semistandard skyline filling, we have
$F(i', j+2)>t$, which, together  with the constraint  that the
entries in row $i'$ are weakly decreasing,
forces  that $F(i', j+2)=t+1$. This verifies \eqref{PPP}.

If column $j+2$ does not contain $t$, then the entry $t+1$ in column $j+2$ is
the pseudo-free entry associated with  $t$. Otherwise, we use the same arguments as above to conclude that
$F(i', j+3)=t+1$. Continuing this procedure, we can eventually  locate
a column index, say $j'$ with $j'>j$, such that $F(i', j')=t+1$ and  column $j'$ of $F$
does not contain $t$, as desired.

We next verify that (1) implies (2). Assume that the   pseudo-free entry
 $t+1$ associated to $t$ is in column $j'$. If $j'=j+1$, then we are done.
 We next consider the case $j'>j+1$. By \eqref{PP-x} in the definition of
 pseudo-free entries, in column $j'-1$, the entry $t+1$ is above the
 entry $t$. The proof
for \eqref{PPP}  implies that $F(i',j'-1)$ must  equal $t+1$. Repeating this procedure
 will lead to $F(i',j+1)=t+1$. So the proof is complete.
\qed

In view of  the proof of Lemma \ref{lem-1},
a typical local configuration of a   pseudo-free entry $t$ and
its associated pseudo-free entry $t+1$ is   illustrated in Figure \ref{Fig-m}, where
the pseudo-free entries are signified in boldface.

\begin{figure}[h]
\begin{center}
\begin{tikzpicture}[scale=0.8]
\draw (1.2,3) node {$t$+$1$}; \draw (2.4,3) node {$t$+$1$};
\draw (3.6,3) node {$t$+$1$}; \draw (4.8,3) node {$\bm{t}$+$\bm{1}$};
\draw (0.2,2.2) node {$\bm{t}$};
\draw (-1.3,2.2) node {row $i$:};
\draw (-1.3,3) node {row $i'$:}; 
\draw (1.2,2.2) node {$t$};
\draw (2.4,1.4) node {$t$}; \draw (3.6,0.8) node {$t$};
\end{tikzpicture}
\end{center}
\vspace{-4mm}
\caption{An illustration of pseudo-free entries. }\label{Fig-m}
\end{figure}

\subsection{The lowering operator}\label{sub-22}

Based on the classification of the entries $t$ and
$t+1$, we can now   define the lowering operator
  $L_{r,t}$ on semistandard
skyline fillings. Note that, according to the flag condition \eqref{ifii},
we naturally have $r\geq t$.

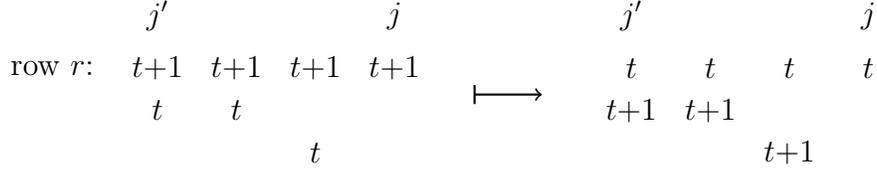
\begin{figure}[h]
\begin{center}
\begin{tikzpicture}[scale=0.7]
\node at (-0.5,3) {row $r$:};
\node at (1.5,4) {$j'$};\node at (6,4) {$j$};
\draw (1.5,3) node {$t$+$1$}; \draw (3,3) node {$t$+$1$};
\draw (4.5,3) node {$t$+$1$}; \draw (6,3) node {$t$+$1$};
\draw (1.5,2.2) node {$t$};
\draw (3,2.2) node {$t$}; \draw (4.5,1.4) node {$t$};
\draw [|->,thick] (7.5,2.5)--(8.8,2.5);
\draw [shift={(9,0)}](1.5,3) node {$t$}; \draw [shift={(9,0)}] (3,3) node {$t$};
\draw [shift={(9,0)}] (4.5,3) node {$t$}; \draw [shift={(9,0)}] (6,3) node {$t$};
\draw (1.5,2.2) [shift={(9,0)}] node {$t$+$1$};
\draw (3,2.2) [shift={(9,0)}] node {$t$+$1$}; \draw [shift={(9,0)}] (4.5,1.4) node {$t$+$1$};

\node at (10.5,4) {$j'$};\node at (15,4) {$j$};
\end{tikzpicture}
\end{center}
\vspace{-5mm}
\caption{An illustration of the lowering operation.}\label{Fig-1}
\end{figure}

Let  $F\in \mathrm{SSF}(\alpha)$,  and denote  $F'=L_{r, t}(F)$.  Consider the
free entries $t+1$ in row $r$ of $F$. If row $r$
of $F$ does not contain any free entry $t+1$, set $F'=F$. Otherwise, locate the
{\it rightmost free entry} $t+1$ in row $r$. Assume  that such a $t+1$ lies in  column $j$.
Locate the smallest column index   $j'<j$, such that for any $j'\leq k<j$,  $F(r, k)=t+1$, and
there is an entry $t$ in  column $k$ lying  below row $r$.
Let   $F'$ be obtained from $F$ by replacing the entry  $t+1$ in column $j$ by $t$, and
then exchanging the  entries $t$ and $t+1$ in   column $k$ for $j'\leq  k<j$. If such column index $j'$ does not exist, then $F'$ is obtained from $F$ by just replacing the entry  $t+1$ in column $j$ by $t$.
This operation is illustrated in Figure \ref{Fig-1}.

\begin{prop}\label{thm-11}
For   $F\in \mathrm{SSF}(\alpha)$, the skyline  filling $L_{r,t}(F)$
also belongs to $\mathrm{SSF}(\alpha)$.
\end{prop}

Clearly, $L_{r,t}(F)$ satisfies items \eqref{ifii} and \eqref{if3i} in the definition of a
semistandard skyline filling. To conclude Proposition \ref{thm-11},
we need to  show that $L_{r,t}(F)$  satisfies items
\eqref{ifi} and \eqref{ifv}, which will be  verified in Lemma
\ref{lem-3} and Lemma \ref{lem-4} respectively.

\begin{lem}\label{lem-3}
For   $F\in \mathrm{SSF}(\alpha)$, the skyline filling $L_{r,t}(F)$ satisfies item \eqref{ifi}.
\end{lem}

\pf
Let $F'=L_{r,t}(F)$,  and let $j, j'$ be the column indices as used in definition of $L_{r,t}$.
For $j'\leq k<j$, assume that the entry $t$ in column $k$ of
$F$ lies in row $r_k$. We need to  show that the entries in row $r$ as well as in row $r_k$ of $F'$ are
weakly decreasing.

Let us first consider row $r$ of $F'$. It suffices to
verify that if $(r,j+1)$ is a box of $D(\alpha)$,
then  $F(r, j+1)\neq t+1$. Suppose otherwise that $F(r, j+1)=t+1$.
Since $F(r,j)=t+1$ and this $t+1$ is free, by definition, the entry $t+1$ in the box $(r, j+1)$ is not pseudo-free.
Recalling that the entry $t+1$ in the box $(r,j)$ is the rightmost  free
entry in row $r$, it follows that
the entry $t+1$ in  $(r, j+1)$ is paired.
Assume that the entry $t$ in column $j+1$ of $F$ is in row $p$.
There are two cases.

Case 1. $p<r$.  Since  row $p$ of $F$ is weakly decreasing and
column $j$ of $F$ does not contain $t$, we have $F(p,j)>t$. Since
column $j$ of $F$ has distinct entries ,  we obtain that $F(p,j)>t+1$,
which, in view of item \eqref{ifv},  leads to $F(p,j+1)>t+1$, contrary to
the assumption that $F(p,j+1)=t$.

Case 2. $p>r$.
By the  arguments
 in first two paragraphs in the  proof of Lemma \ref{lem-1},
we can find a column index $j''>j$ such that  $F(r,j'')=r+1$, column $j''$ of $F$ does not contain
 $t$, and for each $j<k<j''$,
 column $k$ of $F$ contains both $t$ and $t+1$.  Clearly,
the unpaired  entry $t+1$ in column $j''$ is not pseudo-free, so it must be free. This contradicts
the assumption that in row $r$ of $F$,
the entry $t+1$ in column $j$   is  the rightmost    free
entry.

By the above arguments, the assumption that $F(r, j+1)=t+1$ is false, and
hence row $r$ of $F'$ is weakly decreasing.

We next consider  row $r_k$ of $F'$ for $j'\leq k<j$.
By Lemma \ref{BB-f}, we see that
\[r_{j'}\leq r_{j'+1}\leq \cdots \leq r_{j-1}.\]
So we can find column indices   $j'=a_1<a_2<\cdots<a_{m}={j}$   such that for $1\leq h\leq m-1$,
\[r_{a_h}=\cdots=r_{a_{h+1}-1}.\]
Now we need to show that for each $1\leq h\leq m-1$, row $r_{a_h}$ of $F'$ is weakly decreasing.
For    $1<h<m-1$, since column $a_h-1$ of $F$ contains both
$t$ and $t+1$ lying strictly above row $r_{a_h}$, we see that
 $F(r_{a_h}, a_h-1)>t+1$, and so row $r_{a_h}$ of $F'$ is weakly decreasing.

It remains to check that row $r_{j'}=r_{a_1}$ of $F'$ is weakly decreasing.
For simplicity, write $p=r_{j'}$.
We need to verify that $F(p, j'-1)>t$.
Suppose to the contrary that $F(p, j'-1)=t$.
The discussion is divided into two cases.

Case 1. Column $j'-1$ of $F$ does not contain  $t+1$. In this case,
 the entry $t$ in column $j'-1$ and the
 entry $t+1$ in column $j$  would be pseudo-free, leading to a contradiction.

Case 2. Column $j'-1$ of $F$   contains $t+1$. Let $q$ be the row index such that $F(q,j'-1)=t+1$. Since $F(r,j')=t+1$, by Lemma \ref{BB-f}, we see that $q\leq r$. By the choice of  the column index  $j'$,
it follows  that $q\neq r$. So we have    $q<r$.
In light of the proof of  \eqref{PPP},
we obtain that $F(q,j')=t+1$, contrary to the fact $F(r, j')=t+1$.

By the above arguments, the assumption that $F(p, j'-1)=t$ is false.  This completes
the proof.
\qed

\begin{lem}\label{lem-4}
For   $F\in \mathrm{SSF}(\alpha)$, the skyline filling $L_{r,t}(F)$  satisfies
item \eqref{ifv}.
\end{lem}

\pf  Still, let $F'=L_{r,t}(F)$,  and $j, j'$ be the column indices used in the definition of
$L_{r, t}$.
 Let $a'<b'$ be two entries in the same column of $F'$, say column $m$,
such that $a'$ is below $b'$. We need to check that there is an entry  $c'$
in $F'$
to the immediately  right of $b'$ satisfying $a'<c'$.
This is clearly true for
 $m<j'-1$ or $m>j$, since the columns, which are strictly to the right of
 column $j$ or strictly to
 the left of column $j'$, coincide  in $F'$ and $F$. We next
verify the case when $j'-1\leq m\leq j$.

Assume that the entries $a'$ and $b'$ in $F'$ are filled in
the boxes $A$ and $B$, respectively.
Denote by  $a$ and $b$   the entries in $F$ that are filled in
$A$ and  $B$, respectively.

Let us first consider the case $m=j$. Keep in mind that
column $j$ of $F'$ is obtained from column $j$ of $F$
by replacing $t+1$  by $t$.
Since $b'>a'$ in $F'$,  it is easily seen that  $b>a$ in $F$. So
there is an entry $c$  in $F$ immediately
to the right of $b$ such that $a<c$.
Note that column $j+1$  in $F$ and $F'$ coincide.
So the entry   in $F'$ immediately
to the right of $b'$ is also $c$, which is larger than $a'$
by noticing that $a\geq a'$.

We next consider the case $m=j-1$. In this case, column $j-1$ of $F'$ is obtained from
column $j-1$ of $F$
by interchanging  $t+1$  and  $t$. It is readily  checked  that the
only possible pair  that might violate   item \eqref{ifv} would be $a'=t$ and $b'=t+1$.
However, this cannot occur since $t$ lies above $t+1$ in column $j-1$
of $F'$.
Using similar arguments, one can   verify
the cases for $k=j-2, j-3, \ldots, j'-1$. This
completes the proof.
\qed

\subsection{The raising operator}

The raising operator  $R_{r,t}$
is the reverse procedure
of the lowering operator. When implementing the raising operator,
we shall always assume that $r\geq t+1$, since otherwise the resulting
skyline filling would violate the flag condition \eqref{ifii}.

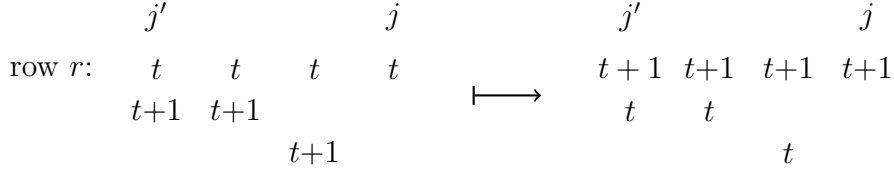
\begin{figure}[h]
\begin{center}
\begin{tikzpicture}[scale=0.7]

\node at (-0.5,3) {row $r$:};
\node at (1.5,4) {$j'$};\node at (6,4) {$j$};

\draw (1.5,3) node {$t$}; \draw (3,3) node {$t$};
\draw (4.5,3) node {$t$}; \draw (6,3) node {$t$};
\draw (1.5,2.2) node {$t$+$1$};
\draw (3,2.2) node {$t$+$1$}; \draw (4.5,1.4) node {$t$+$1$};
\draw [|->,thick] (7.5,2.5)--(8.8,2.5);
\draw [shift={(9,0)}](1.5,3) node {$t+1$}; \draw [shift={(9,0)}] (3,3) node {$t$+$1$};
\draw [shift={(9,0)}] (4.5,3) node {$t$+$1$}; \draw [shift={(9,0)}] (6,3) node {$t$+$1$};
\draw (1.5,2.2) [shift={(9,0)}] node {$t$};
\draw (3,2.2) [shift={(9,0)}] node {$t$}; \draw [shift={(9,0)}] (4.5,1.4) node {$t$};

\node at (10.5,4) {$j'$};\node at (15,4) {$j$};
\end{tikzpicture}
\end{center}
\vspace{-4mm}
\caption{An illustration of the raising operation.}\label{Fig-2}
\end{figure}

Let $F\in \mathrm{SSF}(\alpha)$.
Define  $F'=R_{r, t}(F)$ as follows.
Consider the free entries $t$ in row $r$ of $F$.  If row $r$
of $F$ does not contain any free entry $t$, then set $F'=F$. Otherwise, locate the
{\it leftmost free entry} $t$ in row $r$,
and assume that such a $t$ lies in  column $j$.
Locate the smallest column index   $j'<j$, such that for any $j'\leq k<j$,  $F(r, k)=t$, and
there is an entry $t+1$ in  column $k$ and below row $r$.
Let   $F'$ be obtained from $F$ by replacing the entry  $t$ in column $j$ by $t+1$, and
then exchanging the  entries $t$ and $t+1$ in each column $k$ for $j'\leq  k<j$. If such $j'$ does not exist, then $F'$ is obtained from $F$ by just replacing the entry  $t$ in column $j$ by $t+1$.
This raising operator is illustrated in Figure \ref{Fig-2}.

Using nearly the same arguments as
for Proposition \ref{thm-11},
we can obtain  that the raising operator
is  a map on semistandard skyline fillings.

\begin{prop}\label{thm-22}
For   $F\in \mathrm{SSF}(\alpha)$, the skyline filling $R_{r,t}(F)$
belongs to $\mathrm{SSF}(\alpha)$.
\end{prop}

By the constructions of $L_{r,t}$ and $R_{r,t}$,
we   have the following property.

\begin{coro}\label{coro-x}
Let $F$ be a semistandard skyline filling in $\mathrm{SSF}(\alpha)$.
\begin{itemize}
\item[(1)] If row $r$ of $F$ has a free entry $t+1$, then $R_{r,t}(L_{r,t}(F))=F$;

\item[(2)] If row $r$ of $F$ has a free entry $t$, then $L_{r,t}(R_{r,t}(F))=F$.

\end{itemize}
\end{coro}

\pf
We only give a proof of (1). Let $F\in\mathrm{SSF}(\alpha)$, and $F'=L_{r,t}(F)$. Assume that
 the rightmost free entry $t+1$  in row $r$ of $F$ is in  column $j$.
Then $L_{r,t}$ replaces  this $t+1$  by $t$. By the constructions of
$L_{r,t}$ and $R_{r,t}$, to show that $R_{r,t}(F')=F$,
we need to check that
the entry $t$ in the box $(r,j)$ of $F'$ is still free.
Equivalently, we need to check that the entry $t$ in the box $(r,j)$ of $F'$ is
not pseudo-free.
Suppose to the contrary that it is pseudo-free. By Lemma \ref{lem-1},
there is a row index $r'<r$ such that $F'(r',j+1)=t+1$.
So  $F(r',j+1)=F'(r',j+1)=t+1$, which contradicts   Lemma \ref{BB-f} since
 $F(r,j)=t+1$.  This completes the proof.
\qed

\subsection{The involution $\Phi_{r,t}$ on $\mathrm{SSF}(\alpha)$}

We can
now describe  the  involution $\Phi_{r,t}$, where $r\geq t+1$, on semistandard
skyline fillings. Let $F\in\mathrm{SSF}(\alpha)$.
 Assume that in row $r$ of $F$,
the number of free  entries $t+1$  is  $n_1$, and
the number of free  entries $t$    is  $n_2$.
There are three cases.
\begin{itemize}
\item[(I)] $n_1=n_2$. In this case, set $\Phi_{r, t}(F)=F$.

\item[(II)] $n_1>n_2$. Write $m=n_1-n_2$. Define $\Phi_{r, t}(F)=L_{r,t}^m(F)$ to be the skyline filling
obtained by applying the lowering operator $m$ times to $F$.

\item[(III)] $n_1<n_2$. Write $m'=n_2-n_1$. Define $\Phi_{r, t}(F)=R_{r,t}^{m'}(F)$ to
be the skyline filling
obtained by applying the raising  operator $m'$ times to $F$.
\end{itemize}

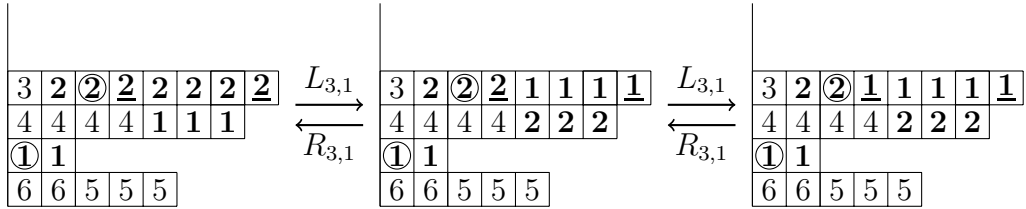
\begin{figure}[h]

\begin{center}
\begin{tikzpicture}[scale=0.45]
\draw (1,-1) -- (1,5);
\draw (1,-1) grid (3,3); \draw (3,-1) grid (6,0);
\draw (3,1) grid (8,3);\draw (7,2) grid (9,3);
\draw (1.5,2.5) node {$3$};
\draw (2.5,2.5) node {$\bm2$}; \draw (3.5,2.5) node {${\bm2}$};
\draw (4.5,2.5) node {$\underline{\bm2}$}; \draw (5.5,2.5) node {${\bm2}$};
\draw (6.5,2.5) node {${\bm2}$}; \draw (7.5,2.5) node {${\bm2}$}; \draw (8.5,2.5) node {$\underline{\bm2}$};
 \draw (1.5,1.5) node {$4$};
\draw (2.5,1.5) node {$4$}; \draw (3.5,1.5) node {$4$};
\draw (4.5,1.5) node {$4$}; \draw (5.5,1.5) node {$\bm1$};
\draw (6.5,1.5) node {$\bm1$};  \draw (7.5,1.5) node {$\bm1$};
 \draw (1.5,0.5) node {$\bm1$};
\draw (2.5,0.5) node {$\bm1$};
 \draw (1.5,-0.5) node {$6$};
\draw (2.5,-0.5) node {$6$}; \draw (3.5,-0.5) node {$5$};
\draw (4.5,-0.5) node {$5$}; \draw (5.5,-0.5) node {$5$};
\draw (1.5,0.5) circle (0.4cm); \draw (3.5,2.5) circle (0.4cm);

\draw [->,thick] (9.5,2)--(11.5,2);
\draw (10.5,2.7) node {$L_{3,1}$};
\draw [<-,thick] (9.5,1.4)--(11.5,1.4);
\draw (10.5,.7) node {$R_{3,1}$};

\draw [shift={(11,0)}] (1,-1) -- (1,5);
\draw[shift={(11,0)}] (1,-1) grid (3,3); \draw [shift={(11,0)}](3,-1) grid (6,0);
\draw[shift={(11,0)}] (3,1) grid (8,3);\draw [shift={(11,0)}](7,2) grid (9,3);
\draw[shift={(11,0)}] (1.5,2.5) node {$3$};
\draw [shift={(11,0)}](2.5,2.5) node {$\bm2$}; \draw [shift={(11,0)}](3.5,2.5) node {${\bm2}$};
\draw [shift={(11,0)}](4.5,2.5) node {$\underline{\bm2}$}; \draw [shift={(11,0)}](5.5,2.5) node {${\bm1}$};
\draw [shift={(11,0)}](6.5,2.5) node {${\bm1}$}; \draw[shift={(11,0)}] (7.5,2.5) node {${\bm1}$}; \draw[shift={(11,0)}] (8.5,2.5) node {$\underline{\bm1}$};
 \draw[shift={(11,0)}] (1.5,1.5) node {$4$};
\draw[shift={(11,0)}] (2.5,1.5) node {$4$}; \draw[shift={(11,0)}] (3.5,1.5) node {$4$};
\draw[shift={(11,0)}] (4.5,1.5) node {$4$}; \draw [shift={(11,0)}](5.5,1.5) node {$\bm2$};
\draw[shift={(11,0)}] (6.5,1.5) node {$\bm2$};  \draw[shift={(11,0)}] (7.5,1.5) node {$\bm2$};
 \draw [shift={(11,0)}](1.5,0.5) node {$\bm1$};
\draw [shift={(11,0)}](2.5,0.5) node {$\bm1$};
 \draw [shift={(11,0)}](1.5,-0.5) node {$6$};
\draw [shift={(11,0)}](2.5,-0.5) node {$6$}; \draw [shift={(11,0)}](3.5,-0.5) node {$5$};
\draw [shift={(11,0)}](4.5,-0.5) node {$5$}; \draw [shift={(11,0)}](5.5,-0.5) node {$5$};
\draw [shift={(11,0)}](1.5,0.5) circle (0.4cm); \draw [shift={(11,0)}](3.5,2.5) circle (0.4cm);

\draw [shift={(11,0)}][->,thick] (9.5,2)--(11.5,2);
\draw [shift={(11,0)}](10.5,2.7) node {$L_{3,1}$};
\draw [shift={(11,0)}] [<-,thick] (9.5,1.4)--(11.5,1.4);
\draw [shift={(11,0)}](10.5,.7) node {$R_{3,1}$};

\draw [shift={(22,0)}] (1,-1) -- (1,5);
\draw[shift={(22,0)}] (1,-1) grid (3,3); \draw [shift={(22,0)}](3,-1) grid (6,0);
\draw[shift={(22,0)}] (3,1) grid (8,3);\draw [shift={(22,0)}](7,2) grid (9,3);
\draw[shift={(22,0)}] (1.5,2.5) node {$3$};
\draw [shift={(22,0)}](2.5,2.5) node {$\bm2$}; \draw [shift={(22,0)}](3.5,2.5) node {${\bm2}$};
\draw [shift={(22,0)}](4.5,2.5) node {$\underline{\bm1}$}; \draw [shift={(22,0)}](5.5,2.5) node {${\bm1}$};
\draw [shift={(22,0)}](6.5,2.5) node {${\bm1}$}; \draw[shift={(22,0)}] (7.5,2.5) node {${\bm1}$}; \draw[shift={(22,0)}] (8.5,2.5) node {$\underline{\bm1}$};
 \draw[shift={(22,0)}] (1.5,1.5) node {$4$};
\draw[shift={(22,0)}] (2.5,1.5) node {$4$};
\draw[shift={(22,0)}] (3.5,1.5) node {$4$};
\draw[shift={(22,0)}] (4.5,1.5) node {$4$};
\draw [shift={(22,0)}](5.5,1.5) node {$\bm2$};
\draw[shift={(22,0)}] (6.5,1.5) node {$\bm2$};  \draw[shift={(22,0)}] (7.5,1.5) node {$\bm2$};
 \draw [shift={(22,0)}](1.5,0.5) node {$\bm1$};
\draw [shift={(22,0)}](2.5,0.5) node {$\bm1$};
 \draw [shift={(22,0)}](1.5,-0.5) node {$6$};
\draw [shift={(22,0)}](2.5,-0.5) node {$6$}; \draw [shift={(22,0)}](3.5,-0.5) node {$5$};
\draw [shift={(22,0)}](4.5,-0.5) node {$5$}; \draw [shift={(22,0)}](5.5,-0.5) node {$5$};
\draw [shift={(22,0)}](1.5,0.5) circle (0.4cm); \draw [shift={(22,0)}](3.5,2.5) circle (0.4cm);

\end{tikzpicture}
\end{center}

\vspace{-4mm}
\caption{An illustration of the involution $\Phi_{r,t}$.}\label{ivl}
\end{figure}

For example, Figure \ref{ivl} illustrates the involution $\Phi_{3,1}$, where $r=3$
and $t=1$. For the leftmost skyline filling, we have $n_1=2, n_2=0$, thus $\Phi_{3,1}=L_{3,1}^2$.

\vspace{2mm}

By Corollary \ref{coro-x}, we  have the following conclusion.

\begin{coro}\label{Con-x}
The map $\Phi_{r, t}$ is an involution on $\mathrm{SSF}(\alpha)$.
Particularly,  $\Phi_{r, t}$ exchanges the numbers of free entries $t$ and $t+1$
in row $r$ of any skyline filling in $\mathrm{SSF}(\alpha)$.
\end{coro}

It can be readily seen that for two distinct row indices $r$ and $r'$,
the columns interfered by    $\Phi_{r, t}$ differ
 from the columns interfered by    $\Phi_{r', t}$, and thus they  are commutative, namely,
\begin{equation}\label{comm}
\Phi_{r', t}\circ \Phi_{r, t}=\Phi_{r, t}\circ \Phi_{r', t}.
\end{equation}

\begin{re}
Consider  a  composition $\alpha=(\alpha_1,\ldots,
\alpha_n)$ whose parts are weakly increasing.
This means its reverse  $\alpha^{\mathrm{rev}}=(\alpha_n,\ldots,
\alpha_1)$ is a partition.
Let $F\in\mathrm{SSF}(\alpha)$. Using item  \eqref{ifv}, 
it is easy to check that
the entries in each  column of $F$ are strictly  increasing.
Let $\overline{F}$ be obtained from $F$
by applying a reflection about the horizontal line. Then each
row of $\overline{F}$ is weakly decreasing, and each column
of $\overline{F}$ is strictly decreasing.
So $\overline{F}$
is a reverse semistandard Young tableau (also called a column-strict plane partition)
of shape $\alpha^{\mathrm{rev}}$ \cite[Chapter 7.10]{Sta}.
In this case, the involution $\Phi_{r,t}$ specifies
to the classical  Bender--Knuth involution \cite{Ben},
which has the form as illustrated in Figure \ref{bk}.
Note that in such a case,   each entry $t$ or $t+1$ is
either paired or free, there are no pseudo-free entries.

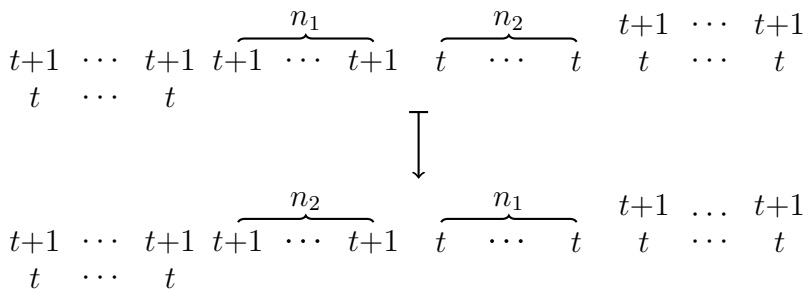
\begin{figure}[h]
\begin{center}
\begin{tikzpicture}[scale=0.6]
\draw (0,2.8) node {$t$+$1$};
\draw (3,2.8) node {$t$+$1$};
\draw (4.5,2.8) node {$t$+$1$};
\draw (7.5,2.8) node {$t$+$1$};
\draw (0,2) node {$t$};
\draw (1.5,2) node {$\cdots$};
\draw (3,2) node {$t$};
\draw (9,2.8) node {$t$};
\draw (12,2.8) node {$t$};
\draw (13.5,2.8) node {$t$};
\draw (16.5,2.8) node {$t$};
\draw (13.5,3.6) node {$t$+$1$};
\draw (15,3.6) node {$\cdots$};
\draw (16.5,3.6) node {$t$+$1$};
\draw[decorate,decoration={brace,raise=4pt},thick] (4.5,3)--(7.5,3);
\draw (6,3.7) node {$n_1$};
\draw[decorate,decoration={brace,raise=4pt},thick] (9,3)--(12,3);
\draw (10.5,3.7) node {$n_2$};
\draw (6,2.8) node {$\cdots$};\draw (1.5,2.8) node {$\cdots$};\draw (15,2.8) node {$\cdots$};\draw (10.5,2.8) node {$\cdots$};

\draw [|->,thick] (8.5,1.7)--(8.5,0.2);

\draw [shift={(0,-4)}](0,2.8) node {$t$+$1$};
\draw [shift={(0,-4)}](3,2.8) node {$t$+$1$};
\draw [shift={(0,-4)}](4.5,2.8) node {$t$+$1$};
\draw [shift={(0,-4)}](6,2.8) node {$\cdots$}; \draw[shift={(0,-4)}] (7.5,2.8) node {$t$+$1$};
\draw [shift={(0,-4)}](0,2) node {$t$};
\draw[shift={(0,-4)}] (1.5,2) node {$\cdots$};
\draw [shift={(0,-4)}](3,2) node {$t$};
\draw [shift={(0,-4)}](9,2.8) node {$t$};
\draw [shift={(0,-4)}](10.5,2.8) node {$\cdots$};
\draw [shift={(0,-4)}](12,2.8) node {$t$};
\draw [shift={(0,-4)}](13.5,2.8) node {$t$};
\draw [shift={(0,-4)}](16.5,2.8) node {$t$};
\draw [shift={(0,-4)}](13.5,3.6) node {$t$+$1$};
\draw[shift={(0,-4)}] (15,3.4) node {$\cdots$};
\draw [shift={(0,-4)}](16.5,3.6) node {$t$+$1$};
\draw[decorate,decoration={brace,raise=4pt},thick] [shift={(0,-4)}](4.5,3)--(7.5,3);
\draw [shift={(0,-4)}](6,3.7) node {$n_2$};
\draw[decorate,decoration={brace,raise=4pt},thick][shift={(0,-4)}] (9,3)--(12,3);
\draw [shift={(0,-4)}](10.5,3.7) node {$n_1$};
\draw [shift={(0,-4)}](6,2.8) node {$\cdots$};
\draw [shift={(0,-4)}](1.5,2.8) node {$\cdots$};
\draw [shift={(0,-4)}](15,2.8) node {$\cdots$};
\draw[shift={(0,-4)}] (10.5,2.8) node {$\cdots$};
\end{tikzpicture}
\end{center}
\vspace{-3mm}
\caption{The Bender--Knuth involution.}\label{bk}
\end{figure}

\end{re}

\section{Application to key polynomials}\label{sec-3}

In this section, we always assume that
$\alpha=(\alpha_1,\alpha_2, \ldots, \alpha_n)\in \mathbb{Z}^n$
is a composition,  $r$ is the {\it first ascent} of $\alpha$
(namely, $\alpha_1\geq \alpha_2\geq \cdots \geq \alpha_r$, and
$\alpha_r<\alpha_{r+1}$), and
$\alpha'$ is obtained
from $\alpha$ by interchanging $\alpha_r$ and $\alpha_{r+1}$.
We aim to apply the involution $\Phi_{r, t}$
 to give a combinatorial interpretation of  the relation
\begin{equation}\label{R-x}
\pi_r\left(
\sum_{F'\in \mathrm{SSF}(\alpha')}x^{F'}\right)=\sum_{F\in \mathrm{SSF}(\alpha)}x^F.
\end{equation}

Using  the involution $\Phi_{r, t}$, we
 define a new involution $\Phi_r$ on  $\mathrm{SSF}(\alpha')$.
For $F'\in\mathrm{SSF}(\alpha')$, let $\Phi_r(F')$ be obtained from $F'$ by
applying  $\Phi_{i,r}$ for each $i\geq r+1$. In view of  \eqref{comm},
this operation does not depend on  the order of the involutions $\Phi_{i,r}$.

By Corollary \ref{Con-x},
we have the following consequence.

\begin{coro}\label{X-1}
For  $F'\in \mathrm{SSF}(\alpha')$, the number of free entries
$r$ (respectively, $r+1$) below row $r$ of
 $F'$ is the same as the number of free entries $r+1$ (respectively, $r$) below
 row $r$ of  $\Phi_r(F')$.
\end{coro}

To prove \eqref{R-x}, we
need three simple lemmas.

\begin{lem}\label{X-2}
Let $F'\in \mathrm{SSF}(\alpha')$. Assume that   row $r$ of $F'$ has $m$ free entries
 $r$. Then
 \begin{equation}\label{E-1}
 x^{F'}+x^{\Phi_r(F')}=x_r^m f(x),
 \end{equation}
 where $f(x)$ is a polynomial symmetric in $x_r$ and $x_{r+1}$.
\end{lem}

\pf
Note that $F'$ and $\Phi_r(F')$ have the same entries other than $r$ and $r+1$.
These entries contribute a common factor to $x^{F'}$ and $x^{\Phi_r(F')}$,
 which does not contain the variables  $x_r$ and $x_{r+1}$.
Note also that  the lowering operator or the raising
operator keeps  paired entries  and  pseudo-free entries unchanged. So
the paired entries and the pseudo-free entries
 together contribute a common factor
of  the form $(x_rx_{r+1})^h$ to $x^{F'}$ and $x^{\Phi_r(F')}$.

Now we consider the free entries in
 $F'$ and $\Phi_r(F')$. First, the entries $r$ and $r+1$ in $F'$ or $\Phi_r(F')$
cannot occur above  row $r$. Second, there are no entries $r+1$ in  row $r$
of $F'$ or $\Phi_r(F')$, and so
the free entries in row $r$ of $F'$ and $\Phi_r(F')$ contribute a common factor
$x_r^m$.
Finally, by Corollary
\ref{X-1},  the total free entries
below row $r$ of
  $F'$ and $\Phi_r(F')$ together generate a term symmetric in  $x_r$ and $x_{r+1}$.
The above analysis allows us to conclude  \eqref{E-1}.
\qed

Let $F'\in \mathrm{SSF}(\alpha')$. Suppose that   row $r$ of $F'$ has $m$ free entries
 $r$. Define $m+1$ semistandard skyline fillings
  $F_0, F_1, \ldots, F_m$ belonging to  $ \mathrm{SSF}(\alpha)$ as follows.
Let $F_0$ be obtained from $F'$ by moving the last $\alpha_{r+1}-\alpha_r$
boxes, together with the filled entries, down to row $r+1$.
It can be
checked  that  $F_0\in\mathrm{SSF}(\alpha)$
by noticing the following two  observations:

\begin{enumerate}[(i)]
\item Let $F\in  \mathrm{SSF}(\alpha)$.
For $1\leq i\leq r$, each box in row $i$ of $T$ is
filled with $i$. The first $\alpha_r$ boxes in row $r+1$ of $F$ are filled
with $r+1$.

\item Let $F'\in  \mathrm{SSF}(\alpha')$.
For $1\leq i<r$ or $i=r+1$, each box in row $i$ of $F'$ is
filled with $i$. The first $\alpha_r$ boxes in row $r$ of $F'$ are filled
with $r$.
\end{enumerate}

For $1\leq k\leq m$, let
\[F_i=R_{r+1, r}^k(F_0)\]
 be the skyline filling obtained  by applying the raising  operator $R_{r+1, r}$ $k$ times to $F_0$.
We call $F_0, F_1,\ldots, F_m$ the {\it derived fillings} of $F'$, and we denote
\[\mathrm{DF}(F')=\{F_0, F_1,\ldots, F_m\}.\]

For example, Figure \ref{df} displays the construction of $\mathrm{DF}(F')$ for a
skyline filling $F'\in\mathrm{SSF}(\alpha')$, where $\alpha=(4,2,6,0,4)$,  $r=2$, and $m=2$.

\begin{figure}[h]
\begin{center}
\begin{tikzpicture}[scale=0.4]
\draw (0,0)--(0,5);
\draw (0,4) grid (4,5);
\draw (0,3) grid (6,4);
\draw (0,2) grid (2,3);
\draw (0,0) grid (4,1);
\draw (0.5,4.5) node {$1$};\draw (1.5,4.5) node {$1$};
\draw (2.5,4.5) node {$1$};\draw (3.5,4.5) node {$1$};
\draw (0.5,3.5) node {$\bm2$};\draw (1.5,3.5) node {$\bm2$};
\draw (2.5,3.5) node {$\bm2$};\draw (3.5,3.5) node {$\bm2$};
\draw (4.5,3.5) node {$\bm2$};\draw (5.5,3.5) node {$1$};
\draw (0.5,2.5) node {$\bm3$};\draw (1.5,2.5) node {$\bm3$};
\draw (0.5,0.5) node {$5$};\draw (1.5,0.5) node {$4$};
\draw (2.5,0.5) node {$4$};\draw (3.5,0.5) node {$\bm3$};
\draw [->,thick] (6.5,2.2)--(8.5,2.2);

\draw [shift={(9,0)}](0,0)--(0,5);
\draw [shift={(9,0)}](0,4) grid (4,5);
\draw [shift={(9,0)}](0,3) grid (2,4);
\draw [shift={(9,0)}](0,2) grid (6,3);
\draw [shift={(9,0)}](0,0) grid (4,1);
\draw [shift={(9,0)}](0.5,4.5) node {$1$};\draw [shift={(9,0)}](1.5,4.5) node {$1$};
\draw [shift={(9,0)}](2.5,4.5) node {$1$};\draw [shift={(9,0)}](3.5,4.5) node {$1$};
\draw [shift={(9,0)}](0.5,2.5) node {$\bm3$};\draw [shift={(9,0)}](1.5,2.5) node {$\bm3$};
\draw [shift={(9,0)}](2.5,2.5) node {$\bm2$};\draw[shift={(9,0)}] (3.5,2.5) node {$\bm2$};
\draw [shift={(9,0)}](4.5,2.5) node {$\bm2$};\draw [shift={(9,0)}](5.5,2.5) node {$1$};
\draw [shift={(9,0)}](0.5,3.5) node {$\bm2$};\draw (1.5,3.5)[shift={(9,0)}] node {$\bm2$};
\draw [shift={(9,0)}](0.5,0.5) node {$5$};\draw [shift={(9,0)}](1.5,0.5) node {$4$};
\draw [shift={(9,0)}](2.5,0.5) node {$4$};\draw [shift={(9,0)}](3.5,0.5) node {$\bm3$};

\draw [shift={(9,0)}][->,thick] (6.5,2.2)--(8.5,2.2);
\draw [shift={(9,0)}]  (7.5,2.8) node {$R_{3,2}$};

\draw [shift={(18,0)}](0,0)--(0,5);
\draw [shift={(18,0)}](0,4) grid (4,5);
\draw [shift={(18,0)}](0,3) grid (2,4);
\draw [shift={(18,0)}](0,2) grid (6,3);
\draw [shift={(18,0)}](0,0) grid (4,1);
\draw [shift={(18,0)}](0.5,4.5) node {$1$};\draw [shift={(18,0)}](1.5,4.5) node {$1$};
\draw [shift={(18,0)}](2.5,4.5) node {$1$};\draw [shift={(18,0)}](3.5,4.5) node {$1$};
\draw [shift={(18,0)}](0.5,2.5) node {$\bm3$};\draw [shift={(18,0)}](1.5,2.5) node {$\bm3$};
\draw [shift={(18,0)}](2.5,2.5) node {$\bm3$};\draw[shift={(18,0)}] (3.5,2.5) node {$\bm2$};
\draw [shift={(18,0)}](4.5,2.5) node {$\bm2$};\draw [shift={(18,0)}](5.5,2.5) node {$1$};
\draw [shift={(18,0)}](0.5,3.5) node {$\bm2$};\draw (1.5,3.5)[shift={(18,0)}] node {$\bm2$};
\draw [shift={(18,0)}](0.5,0.5) node {$5$};\draw [shift={(18,0)}](1.5,0.5) node {$4$};
\draw [shift={(18,0)}](2.5,0.5) node {$4$};\draw [shift={(18,0)}](3.5,0.5) node {$\bm3$};

\draw [shift={(18,0)}][->,thick] (6.5,2.2)--(8.5,2.2);
\draw [shift={(18,0)}]  (7.5,2.8) node {$R_{3,2}$};

\draw [shift={(27,0)}](0,0)--(0,5);
\draw [shift={(27,0)}](0,4) grid (4,5);
\draw [shift={(27,0)}](0,3) grid (2,4);
\draw [shift={(27,0)}](0,2) grid (6,3);
\draw [shift={(27,0)}](0,0) grid (4,1);
\draw [shift={(27,0)}](0.5,4.5) node {$1$};\draw [shift={(27,0)}](1.5,4.5) node {$1$};
\draw [shift={(27,0)}](2.5,4.5) node {$1$};\draw [shift={(27,0)}](3.5,4.5) node {$1$};
\draw [shift={(27,0)}](0.5,2.5) node {$\bm3$};\draw [shift={(27,0)}](1.5,2.5) node {$\bm3$};
\draw [shift={(27,0)}](2.5,2.5) node {$\bm3$};\draw[shift={(27,0)}] (3.5,2.5) node {$\bm3$};
\draw [shift={(27,0)}](4.5,2.5) node {$\bm3$};\draw [shift={(27,0)}](5.5,2.5) node {$1$};
\draw [shift={(27,0)}](0.5,3.5) node {$\bm2$};\draw (1.5,3.5)[shift={(27,0)}] node {$\bm2$};
\draw [shift={(27,0)}](0.5,0.5) node {$5$};\draw [shift={(27,0)}](1.5,0.5) node {$4$};
\draw [shift={(27,0)}](2.5,0.5) node {$4$};\draw [shift={(27,0)}](3.5,0.5) node {$\bm2$};
\draw (2,-1.2) node {$F'$}; \draw (11,-1.2) node {$F_0$}; \draw (20,-1.2) node {$F_1$}; \draw (29,-1.2) node {$F_2$};
\end{tikzpicture}
\end{center}
\vspace{-4mm}
\caption{An illustration of the construction of  $\mathrm{DF}(F')$.}\label{df}
\end{figure}

\begin{lem}\label{X-3}
For $F'\in \mathrm{SSF}(\alpha')$, we have
\begin{equation}\label{E-2}
\pi_r\left(x^{F'}+x^{\Phi_r(F')}\right)=\sum_{F\in \mathrm{DF}(F')} x^{F}+
\sum_{F\in   \mathrm{DF}(\Phi_r(F'))} x^{F}.
\end{equation}
\end{lem}

\pf For a polynomial  $f_1(x)$
symmetric in $x_r$ and $x_{r+1}$,
it is easy to check that for any polynomial $f_2(x)$,
\[\partial_r(f_1(x)\, f_2(x))=f_1(x)\,\partial_r(f_2(x)).\]
So, by Lemma \ref{X-2},
\begin{align}
\pi_r\left(x^{F'}+x^{\Phi_r(F')}\right)&=\partial_r(x^{m+1}f(x))=\partial_r(x^{m+1})\,f(x)\nonumber\\[5pt]
&=(x_r^m+x_r^{m-1}x_{r+1}+\cdots+x_{r+1}^m)f(x).\label{E-4}
\end{align}
In view of  the construction of the derived fillings of $F'$, it is easily
seen that \eqref{E-4} coincides with the  right-hand side of \eqref{E-2}.
\qed

\begin{lem}\label{X-8}
The set $\mathrm{SSF}(\alpha)$ is a disjoint union of $\mathrm{DF}(F')$,
where $F'$ runs over semistandard skyline fillings in $\mathrm{SSF}(\alpha')$. That is,
\[\mathrm{SSF}(\alpha)=\biguplus_{F'\in\mathrm{SSF}(\alpha')}\mathrm{DF}(F').\]
\end{lem}

\pf
It is easy to see that the sets  $\mathrm{DF}(F')$ are disjoint.
On the other hand, given a skyline filling $F\in \mathrm{SSF}(\alpha)$,
the corresponding skyline filling  $F'\in \mathrm{SSF}(\alpha')$ such that
$F\in \mathrm{DF}(F')$ can be  constructed as follows.
Suppose that there are $k$ free entries $r+1$ in row $r$
of $F$. Let $\overline{F}$ be obtained from  $F$ by applying the lowering operator $L_{r+1,r}$
 $k$ times. Let $F'$ be obtained from $\overline{F}$ by moving the last $\alpha_{r+1}-\alpha_{r}$
boxes, together with the filled entries, up to row $r$. It is
routine to check that $F'$ belongs to $\mathrm{SSF}(\alpha')$.
Moreover, it is easy to see that $F$ belongs to $\mathrm{DF}(F')$.
This completes the proof.
\qed

Lemma \ref{X-3}
provides an   algorithm to generate inductively the
semistandard skyline fillings for key polynomials.
By  Lemma \ref{X-3}, we obtain that
\begin{equation}\label{R-0}
\pi_r\left(
\sum_{F'\in \mathrm{SSF}(\alpha')}x^{F'}\right)=
\sum_{F'\in \mathrm{SSF}(\alpha')} \sum_{F\in \mathrm{DF}(F')}x^{F}.
\end{equation}
By Lemma \ref{X-8}, the right-hand side of
\eqref{R-0} is the same as the right-hand side of \eqref{R-x}.
This establishes  relation \eqref{R-x}.

\vspace{.2cm} \noindent{\bf Acknowledgments.}
We are grateful to Oliver Pechenik and Dominic Searles for explaining the
connections among several combinatorial structures for the construction
of  key polynomials.
This work was
supported by the National Science Foundation of China (Grant No. 11971250, 12071320)
and  the Sichuan Science and Technology Program (Grant No. 2020YJ0006).

\footnotesize{

\textsc{Neil J.Y. Fan, Department of Mathematics, Sichuan University, Chengdu 610064, P.R. China.}
Email address: fan@scu.edu.cn

\textsc{Peter L. Guo, Center for Combinatorics, Nankai University, Tianjin 300071, P.R. China.}
Email address: lguo@nankai.edu.cn

\textsc{Nicolas Y. Liu, Center for Combinatorics, Nankai University, Tianjin 300071, P.R. China.}
Email address: yiliu@mail.nankai.edu.cn

}

\end{document}